\documentclass[12pt]{amsart}
\usepackage{amsxtra,latexsym,amssymb}
\usepackage{epsfig,rotate,amsthm}
\usepackage{hyperref}

\def\qed{{\hfill $\Box$}}

\def\N{{\mathbb N}}
\def\Z{{\mathbb Z}}
\def\C{{\mathbb C}}
\def\Q{{\mathbb Q}}

\def\Nil{U_{r,s}^{+}(\mathfrak g)}
\def\U{{U_{r,s}^{+}({\mathfrak s}\mathfrak{l}_{3})}}
\def\V{U_{r, s}^{+}(B_{2})}
\theoremstyle{theorem}
\newtheorem{thm}{Theorem}[section]
\newtheorem{cor}{Corollary}[section]
\newtheorem{prop}{Proposition}[section]

\newtheorem{lem}{Lemma}[section]
\theoremstyle{definition}
\newtheorem{defn}{Definition}[section]
\theoremstyle{remark}

\begin{document}
\title[ Prime Ideals and Automorphims of $\V$ ]{The Prime ideal Stratification and The Automorphism Group of $U^{+}_{r,s}(B_{2})$}
\author[X. Tang]{Xin Tang}
\address{Department of Mathematics \& Computer Science\\
Fayetteville State University\\
1200 Murchison Road, Fayetteville, NC 28301}
\email{xtang@uncfsu.edu}        
\thanks{This research project is partially supported by the ISAS mini-grant at Fayetteville State University}
\date{\today}
\begin{abstract}
Let ${\mathfrak g}$ be a finite dimensional complex simple Lie algebra, and let $r,s\in \C^{\ast}$ be transcendental over $\Q$ such that $r^{m}s^{n}=1$ implies $m=n=0$. We will obtain some basic properties of the two-parameter quantized enveloping algebra $\Nil$. In particular, we will verify that the algebra $\Nil$ satisfies many nice properties such as having normal separation, catenarity and Dixmier-Moeglin equivalence. We shall study a concrete example, the algebra $\V$ in detail. We will first determine the normal elements, prime ideals and primitive ideals for the algebra $\V$, and study their stratifications. Then we will prove that the algebra automorphism group of the algebra $\V$ is isomorphic to $(\C^{\ast})^{2}$. 
\end{abstract}
\maketitle       

\section*{Introduction}

The two-parameter quantized enveloping algebra $U_{r,s}({\mathfrak s}\mathfrak{l}_{n})$ has been studied in \cite{BW1, BW2, Tak}. The two-parameter quantized enveloping algebras $U_{r,s}(\mathfrak g)$ have been studied for the finite dimensional complex simple Lie algebras $\mathfrak g$ of other types in the literatures \cite{BGH, HP} and the references therein. It is easy to see that two-parameter quantized enveloping algebras $U_{r,s}(\mathfrak g)$ are close analogues of the one-parameter quantized enveloping algebras $U_{q}(\mathfrak g)$. On the one hand, they share many common features with their one-parameter analogues. For instance, the two-parameter quantized enveloping algebras $U_{r,s}(\mathfrak g)$ are also Hopf algebras with natural Hopf algebra structures, and they admit triangular decompositions. Furthermore, these Hopf algebras can be realized as the Drinfeld doubles of their Hopf sub-algebras. On the other hand, they do have some different features in both the structure and representation theory \cite{BGH, BKL, BW1, BW2}. 

In order to study the structure of $U_{r,s} (\mathfrak g)$, one is led to first study the two-parameter quantized enveloping algebras $U_{r,s}^{+}(\mathfrak g)$. In references \cite{Rei, T1}, these algebras have been investigated from the point of view of two-parameter Ringel-Hall algebras. Indeed, the algebra $\Nil$ can be presented as an iterated skew polynomial ring, and 
a PBW basis has been constructed for $\Nil$. Furthermore, all the prime ideals of the algebra $\Nil$ are proved to be completely prime ideals based on a mild condition on the parameters $r,s$. 

Historically, there has been much interest in the study of prime ideals and automorphism 
group of quantum algebras including the quantized enveloping algebra $U_{q}(\mathfrak g)$ and their subalgebras, for example \cite{AC, AD, AD1, M}.  Recently, the automorphism group 
problem has been settled for the algebra $U_{q}^{+}(\mathfrak g)$ in some special cases. 
The automorphism group of $U_{q}^{+}({\mathfrak s}\mathfrak{l}_{3})$ was settled by 
Alev and Dumas in \cite{AD} and by Caldero independently in \cite{C2}. Furthermore, the automorphism group of the algebra $U_{q}^{+}(\mathfrak{s}\mathfrak{l}_{4})$ was 
determined in \cite{LL}; and the automorphism group of the algebra $U_{q}^{+}(B_{2})$ 
was determined in \cite{Launois}. In the latter two cases, the central elements and the normal elements (or equivalently the prime and primitive ideals) have played an important role in the determination of the automorphism group. Unfortunately, it still remains a difficult question to determine the automorphism group for any general $U_{q}^{+}(\mathfrak g)$, although it was conjectured in \cite{AD1} that the automorphism group is isomorphic to the semi-direct product of a torus with the group of algebra automorphisms induced by the symmetries of the Dynkin diagram. To get a better picture on the two-parameter quantized enveloping algebra $\Nil$, it is natural to study the prime ideals and the algebra automorphisms of the algebra $\Nil$. Once again, the determination of the normal elements might contribute to a better understanding of the prime and primitive ideals of $\Nil$, and thus 
the automorphism group of $\Nil$.

In this paper, as one step forward, we will first derive some background information on 
prime and primitive ideals of the algebra $\Nil$. Since the defining relations of the algebra 
$\Nil$ are homogeneous, there is a torus (denoted by $\mathcal{H}$) acting on the algebra 
$\Nil$. Therefore, we can apply Goodearl-Letzter stratification theory \cite{Goodearl, GLet} 
to study the prime and primitive spectra of the algebra $\Nil$. Indeed, we will prove that the number of $\mathcal{H}-$invariant prime ideals of the algebra $\Nil$ is equal to the cardinality 
of the Weyl group associated to the Lie algebra $\mathfrak g$. We will use the fact that the algebra $\Nil$ is indeed a cocycle twist of the one-parameter quantized enveloping algebra $U_{q}^{+}(\mathfrak g)$, whose $\mathcal{H}-$invariant prime ideals were proved by 
Gorelik \cite{Gorelik} to be parameterized by the elements of the Weyl group. We will also 
verify that the algebra $\Nil$ have nice properties such as normal separation and being 
catenary. Additionally, we will show the Dixmier-Moeglin equivalence holds for the 
algebra $\Nil$; thus, the primitive ideals of $\Nil$ can be recognized as the maximal 
elements  in each prime stratum. 

We will study a particular case, the algebra $\V$ in terms of its normal elements, prime ideals, primitive ideals, and automorphism group. We will first determine all the normal elements of $\V$. As a result, we determine all the prime and primitive ideals of $\V$ and describe their stratifications. In particular, we will prove that the automorphism group of the algebra $\V$ is a torus of rank two. 

As in the one-parameter case, one may conjecture that the automorphism group of $U^{+}_{r,s}(\mathfrak g)$ is isomorphic to the semi-direct product 
of a torus with the group of algebra automorphisms induced by the symmetries of the corresponding Dynkin diagram. Both the result in \cite{CL} and the result in the current paper serve as positive evidences toward the truth of this conjecture. In a paper \cite{T4} which is currently under preparation, we will address the automorphisms and derivations of the two-parameter quantized enveloping algebra $U_{r,s}^{+}({\mathfrak s}\mathfrak{l}_{4})$. In particular, we shall manage to prove that such a conjecture is true for the algebra $U_{r,s}^{+}({\mathfrak s}\mathfrak{l}_{4})$.

The paper is organized as follows. In section 1, we recall the definition and some basic facts about the two-parameter quantized enveloping algebra $U_{r,s}^{+}(\mathfrak g)$. In Section 2, we will study the normal elements, prime ideals, and primitive ideals for $\Nil$. In Section 3, we treat a particular case, the algebra $\V$.

\section{The two-parameter quantized enveloping algebra $U_{r,s}^{+}(\mathfrak g)$}

In this section, we recall the definition of the two-parameter quantized enveloping algebra $\Nil$ and its Ringel-Hall algebra realization. In particular, we will recall how the algebra $\Nil$ is presented as an iterated skew polynomial ring and the construction of its PBW basis. 
\subsection{The definition of the algebra $\Nil$}
Let $\mathfrak g$ be a finite dimensional complex simple Lie algebra. The two-parameter quantum group $U_{r,s}(\mathfrak g)$ associated to the Lie algebra $\mathfrak g$ has been studied by many authors in the literatures \cite{BW1, BW2, BGH} therein. Moreover, the subalgebras $U_{r,s}^{\geq 0}(\mathfrak g)$ and $\Nil$ of $U_{r,s}(\mathfrak g)$ have been further studied in \cite{Rei, T1} from the point of view of two-parameter Ringel-Hall algebras. In this subsection, we will recall the definition of $\Nil$ and some of its basic properties. 

Let $n$ be the rank of the Lie algebra $\mathfrak g$ and let $I=\{1, 2, \cdots, n\}$. Let us denote by $C=(a_{ij})_{i,j\in I}$ the Cartan matrix associated to the Lie algebra $\mathfrak g$. Let $\{d_{i}\mid i\in I\}$ be a set of relatively prime positive integers such that $d_{i}a_{ij}=d_{j}a_{ji}$ for $i,j\in I$. We will choose two complex numbers $r, s\in \C^{\ast}$ which are transcendental over $\Q$ in such a way that $r^{m}s^{n}=1$ implies that  $m=n=0$. And we will always set $r_{i}=r^{d_{i}}, s_{i}=s^{d_{i}}$.

Let us denote by $\langle-,-\rangle$ the Ringel form (or Euler form) defined on the root lattice $\mathcal{Q}\cong \Z^{n}$ where $n$ is the rank of the Lie algebra $\mathfrak g$. Recall that this bilinear form is defined as follows
\[
\langle i,j  \rangle \colon =\langle \alpha_{i},\alpha_{j} \rangle
=\left
\{\begin{array}{cc}d_{i}a_{ij}, & {\rm if}\, i<j,\\
d_{i},& {\rm if} \,i=j,\\
0, & {\rm if } \,i>j.
\end{array}
\right.
\]

Note that there is a uniform definition proposed for the two-parameter quantum group $U_{r,s}(\mathfrak g)$ in \cite{HP}, and we will recall the following definition of the two-parameter quantized enveloping algebras $U_{r,s}^{+}(\mathfrak g)$ using the notation in \cite{HP}.

\begin{defn} (See also \cite{BGH, BW1}) The two-parameter quantized enveloping algebra $U_{r,s}^{+}(\mathfrak g)$ is defined to be the $\C-$algebras generated by the generators $e_{i}$ subject to the following relations:
\begin{eqnarray*}
\sum_{k=0}^{1-a_{ij}}(-1)^{k}\binom{1-a_{ij}}{k}_{r_{i}s_{i}^{-1}}c_{ij}^{(k)}
e_{i}^{1-a_{ij}-k}e_{j}e_{i}^{k}=0,\quad
(i\neq j)
\end{eqnarray*}
where $c_{ij}^{(k)}=(r_{i}s_{i}^{-1})^{\frac{k(k-1)}{2}}r^{k\langle
  j,i\rangle}s^{-k\langle i,j\rangle }$
for $i\neq j$, and for a symbol $v$, we set up the following notation:
\begin{eqnarray*}
(n)_{v}=\frac{v^{n}-1}{v-1},\quad
(n)_{v}!=(1)_{v}(2)_{v}\cdots(n)_{v},\\
\binom{n}{k}_{v}=\frac{(n)_{v}!}{(k)_{v}!(n-k)_{v}!},\quad \text{for}\, n\geq
k\geq 0,
\end{eqnarray*}
and $(0)_{v}!=1$.
\end{defn}

In the case of the Lie algebra $\mathfrak g =sl_{n+1}$, the two-parameter quantized enveloping algebra $U_{r,s}^{+}({\mathfrak s}\mathfrak{l}_{n+1})$ is indeed generated by the generators $e_{1}, e_{2}, \cdots, e_{n}$ subject to the following relations:
\begin{eqnarray*}
e_{i}^{2}e_{i+1}-(r+s)e_{i}e_{i+1}e_{i}+rse_{i+1}e_{i}^{2}=0,\\
e_{i}e_{i+1}^{2}-(r+s)e_{i+1}e_{i}e_{i+1}+rse_{i+1}^{2}e_{i}=0
\end{eqnarray*}
for $i=1, 2,\cdots, n-1$. It is obvious that the two-parameter quantized enveloping algebra $U_{r,s}^{+}(\mathfrak g)$ is actually a subalgebra of the two-parameter quantum group $U_{r,s}(\mathfrak g)$.

For later on purpose, we need to introduce a gradation on the algebra $\Nil$. Let us denote by $\Z^{n}$ the free abelian group of rank $n$, with a basis denoted by $z_{1}, z_{2}, \cdots, z_{n}$. Given any element ${\bf a} \in \Z^{n}$, say ${\bf a}=\sum a_{i}z_{i}$, we shall set $|{\bf a}|=\sum a_{i}$.  By assigning to the generator $e_{i}$  the degree $z_{i}$, the algebra $U^{+}_{r,s}(\mathfrak g)$ shall become a $\Z^{n}-$graded algebra. Given any ${\bf a}\in \Z^{n}$, we shall denote by $U_{r,s}^{\pm}(\mathfrak g)_{\bf a}$ the set of homogeneous elements of degree ${\bf a}$ in $U^{\pm}_{r,s}(\mathfrak g)$. It is easy to see that we have the following decomposition of the algebra $\Nil$:
\[
U_{r,s}^{+}(\mathfrak g)=\bigoplus_{\bf a}U_{r,s}^{+}(\mathfrak g)_{\bf a}.
\] 

\subsection{Two-parameter Ringel-Hall algebra $H_{r,s}(\Lambda)$}
Let $k$ be a finite field. It is well-known that there is a finite dimensional Hereditary $k-$algebra $\Lambda$ associated to the Lie algebra $\mathfrak g$. We denote by $\mathcal{P}$ the set of all isomorphism classes of finite dimensional $\Lambda-$modules. Thanks to the existence of Hall polynomials, a two-parameter twisted generic Ringel-Hall algebra $H_{r,s}(\Lambda)$ has been constructed in \cite{Rei, T1}. The algebra $H_{r,s}(\Lambda)$ has been successfully used to study the algebra $\Nil$. More general work along the multiparameter quantized enveloping algebras and their realizations can be found in \cite{Pei}. First of all, we will recall the construction  of $H_{r,s}(\Lambda)$, and some of its basic properties here. 

Let us denote by $H_{r,s}(\Lambda)$ the $\C-$linear space spanned by the set $\{u_{\alpha}\mid \alpha \in \mathcal{P}\}$, and let us define a multiplication on the $\C-$linear space $H_{r,s}(\Lambda)$ as follows:
\[
u_{\alpha} u_{\beta} =\sum_{\lambda \in \mathcal{P}}
s^{-\langle \alpha ,\beta \rangle}F^{u_{\lambda} }_{u_{\alpha} u_{\beta}}(rs^{-1})u_{\lambda},\quad \text{for
  any}\, \alpha, \beta \in \mathcal{P}.
\]

Then it is easy to see that $H_{r,s}(\Lambda)$ is an associative $\C-$algebra under the above multiplication.  In addition, $H_{r,s}(\Lambda)$ is a graded algebra whose elements are graded by the dimension vectors. According to \cite{T1}, we can further choose elements $X_{i}$ in $H_{r,s}(\Lambda)$ for $i=1, 2, \cdots, m$, where $m$ is the number of positive roots, such that the following result is true.
\begin{thm}
[Theorem 2.3.2 in \cite{T1}] The monomials $X_{1}^{\alpha(1)}\cdots X_{m}^{
  \alpha(m)}$ with $\alpha(1), \cdots, \alpha(m)\in \N_{0}$ form a
  $\C-$basis of the algebra $H_{r,s}(\Lambda)$; and for $i<j$, we have  
\begin{eqnarray*}
X_{j}X_{i}&= & r^{\langle \dim X_{i}, \dim X_{j}\rangle }s^{-\langle \dim
  X_{j}, \dim X_{i}\rangle } X_{i}X_{j}\\
&& + \sum _{I(i, j)} c(a_{i+1},
  \cdots, a_{j-1}) X_{i+1}^{a_{i+1}}\cdots
  X_{j-1}^{a_{j-1}}
\end{eqnarray*}
with coefficients $c(a_{i+1},\cdots,a_{j-1})$ in $\Q(r,s)$. Here the 
index set $I(i,j)$ is the  set of sequences $(a_{i+1}, \cdots
a_{j-1})$ of natural numbers such that $\sum_{t=i+1}^{j-1}a_{t}{\bf a_{t}=a_{i}+a_{j}}$.
\end{thm}
\qed

Note that by a prime ideal of $A$, we mean a proper two-sided ideal $I\subset A$ such that $aAb \subset I$ implies $a\in I$ or $b\in I$. A prime ideal $P$ is called completely prime if $A/P$ is a domain, or equivalently $ab\in I$ implies $a\in P$ or $b\in P$. The previous theorem implies that $H_{r,s}(\Lambda)$ can be presented as an iterated skew polynomial ring, and thus it has a PBW basis. Thanks to \cite{GLet}, we know that all prime ideals of $H_{r,s}(\Lambda)$ are completely prime under the condition that the multiplicative group generated by $r,s$ is torsion-free. Thanks 
to the results in \cite{Rei,T1}, one further knows that the two-parameter quantized enveloping algebra $\Nil$ is indeed isomorphic to the two-parameter Ringel-Hall algebra $H_{r,s}(\Lambda)$. We should also mention that this isomorphism is indeed a graded isomorphism. Therefore, we shall have the following result.
\begin{thm}
The algebra $\Nil$ is an iterated skew polynomial ring. In particular, all prime ideals of $\Nil$ are completely prime.  
\end{thm}
\qed

In addition, one should note that the two-parameter quantized enveloping algebra $\Nil$ is indeed a cocycle twist of the one-parameter quantized enveloping algebra $U_{q}^{+}(\mathfrak g)$. Let us explain in a little of detail. Note that the one-parameter quantized enveloping algebra $U_{q}^{+}(\mathfrak g)$ is isomorphic to the one-parameter Ringel-Hall algebra $H_{v,v^{-1}}(\Lambda)$ with $v^{2}=q$; and the two-parameter quantized enveloping algebra $\Nil$ is isomorphic to the two-parameter Ringel-Hall algebra $H_{r,s}(\lambda)$. Both isomorphisms are indeed graded isomorphisms. In addition, it is obvious to see that both the algebra $H_{v,v^{-1}}(\Lambda)$ and the algebra $H_{r,s}(\Lambda)$ are graded by the same group, and the algebra $H_{r,s}(\Lambda)$ is a cocycle twist of the $H_{v,v^{-1}}(\Lambda)$.  Indeed, one can 
define the cocycle twist of the algebra $H_{v,v^{-1}}(\Lambda)$ as follows
\[
x \ast y = (v^{-1}s^{-1})^{\langle \mu, \nu \rangle} xy
\]
for any $x\in \Nil_{\mu},  y\in \Nil_{\nu}$.  Then it is easy to see that this cocycle twist (the algebra with the new multiplication $\ast$) is isomorphic to the algebra 
$H_{r,s}(\Lambda)$ under the condition that $rs^{-1}=v^{2}$.

\subsection{The normal elements of the algebra $\Nil$}
Recall the central elements and normal elements of the one-parameter quantized enveloping algebra $U_{q}^{+}(\mathfrak g)$ were described by Caldero in \cite{C1,C2}. And the central elements and normal elements have played an important role in the study of prime ideals and automorphisms of the algebra $U_{q}^{+}(\mathfrak g)$. It is expected that the central elements 
and normal elements of $\Nil$ shall play an equally important role in the determination of prime ideals  and automorphisms of the algebra $\Nil$. 

In this subsection, we shall derive some preliminary information on the normal elements of the algebra $\Nil$. In particular, we will prove that the normal elements of $\Nil$ can be described as $r-s-$central elements of $\Nil$. We need to adopt an approach used by Caldero \cite{C2} for the one-parameter quantized enveloping algebra $U_{q}^{+}(\mathfrak g)$. The following definition mimics the one stated in \cite{C2}.
\begin{defn}
Let $A$ be a $\C-$algebra and $r,s\in \C^{\ast}$ be transcendental over $\Q$ such that $r^{m}s^{n}=1$ implies $m=n=0$. Two elements $u,v\in A $ are called $(r,s)-$commuting if we have $uv=r^{m}s^{n}vu$ for some integers $m,n$. An element in $A$ is called $(r,s)-$central if it is $(r,s)-$commuting with the natural generators of $A$.
\end{defn}

First of all, we have the following lemma.
\begin{lem}
The normal elements of the algebra $\Nil$ are exactly the $(r,s)-$central elements.
\end{lem}
{\bf Proof:} (The proof is an adoption of the one used by Caldero in \cite{C2}). First of all, it is obvious that the $r-s-$central elements of $\Nil$ are normal elements of $\Nil$. Now let $a$ be a non-zero normal element of $\Nil$. Suppose that $a=\sum_{\beta\in \mathcal{Q}^{+}} a_{\beta}$ is a decomposition of the element $a$ in the algebra $\Nil=\bigoplus_{\beta\in \mathcal{Q}^{+}}\Nil_{\beta}$ where $\mathcal{Q}^{+}$ is the positive part of the root lattice. 
Since $a$ is a normal element of $\Nil$, we have $aX_{i}=X_{i}^{\prime} a$, where $X_{i}$ are root vectors corresponding to the simple roots of the Lie algebra $\mathfrak g$. 
By considering the weights, we shall have $a_{\beta}X_{i}=X_{i}^{\prime} a_{\beta}$ for some $X_{i}^{\prime}$. Therefore, we shall have $X_{i}^{\prime}=
\lambda_{i}X_{i}$ for some $\lambda_{i}\in \C$. In particular, we shall have $\lambda_{i}=r^{m_{i}}s^{n_{i}}$. Therefore,  we have finished the proof of the statement.
\qed

\subsection{Quantum unique factorization domains}
The notion of a noncommutative unique factorization domain was introduced by Chatters \cite{Chatters} and were further studied in \cite{CJ}. Recall that an element $p$ of a 
noetherian domain $R$ is said to be prime if $pR = Rp$ and $pR$ is a height one prime ideal of $R$ and $R/pR$ is an integral domain. A noetherian domain $R$ is called a unique 
factorization domain (noetherian UFD) provided that $R$ has at least one height-one prime ideal, and every height-one prime ideal is generated by a prime element. 

Many examples including universal enveloping algebras of finite dimensional solvable 
Lie algebras and their one-parameter quantizations are noetherian UFDs. As a matter of fact, 
it has been proved that many iterated skew polynomial rings are UFDs under some mild conditions. 
In particular, the result ({\bf Theorem 3.7} in \cite{LLR}) states that if $A$ is a torsion-free CGL-extension (Cauchon-Gooderal-Letzter extension) \cite{GLet, Cauchon}, 
then $A$ is a UFD. To proceed, we recall the following definition of the CGL extensions from reference \cite{LLR}.
\begin{defn}
[Definition 3.1 in \cite{LLR}] An iterated skew polynomial extension
$A = k[x_{1} ][x_{2} ; \sigma_{2}, \delta_{2} ] \cdots  [x_{n} ; \sigma_{n}, \delta_{n}]$
is said to be a CGL extension (after Cauchon, Goodearl and Letzter) provided that the
following list of conditions are satisfied:
\begin{itemize}
\item With $A_{j} := k[x_{1} ][x_{2}; \sigma_{2}, \delta_{2} ] \cdots [x_{j}; \sigma_{j} , \delta_{ j}]$ for each $1 \leq j\leq n$, each $\sigma_{j}$ is a $k-$automorphism of $A_{j−1}$, each $\delta_{j}$  is a locally nilpotent $k-$linear $\sigma_{j}-$derivation of $A_{j−1}$ , and there exist
nonroots of unity $q_{j} \in k^{\ast}$ with $\sigma_{j}\delta_{j}=q_{j}\delta_{j}\sigma_{j}$;

\item For each $i < j$ there exists a $\lambda_{ji}$ such that $\sigma_{j}(x_{i})=\lambda_{ji}x_{i}$;

\item There is a torus $\mathcal{H} = (k ^{\ast})^{r}$ acting rationally on $A$ by $k-$algebra automorphisms;

\item The $x_{i}$ for $1\leq i\leq n$ are $\mathcal{H}-$eigenvectors;

\item There exist elements $h_{1}, \cdots, h_{n} \in \mathcal{H}$  such that $h_{j} (x_{i})=\sigma_{j} (x_{i} )$ for $j > i$ and such that the $h_{j} -$eigenvalue of $x_{j}$  is not a root of unity.
\end{itemize}
\end{defn}

Now it is easy to verify that we have the following result about the algebra $\Nil$.
\begin{thm}
The algebra $\Nil$ is a torsion-free $CGL-$extension. In particular, the algebra $\Nil$ is a noetherian UFD.
\end{thm}
\qed

\section{Prime and Primitive ideals of the algebra $\Nil$ }

\subsection{The $\mathcal{H}$-stratification of $Spec(\Nil)$ and $Prim(\Nil)$}
In this subsection, we will apply the stratification theory to the prime spectrum  and primitive spectrum of the algebra $\Nil$. We denote  by $Spec(\Nil)$ the set of all prime ideals of the algebra $\Nil$. And we denote  by $Prim(\Nil)$ the set of all primitive ideals of the algebra $\Nil$. We first recall some results on the stratification theory developed by Goodearl and Letzter \cite{GLet, Goodearl}. As a matter of fact, the stratification theory produces partitions of the prime and primitive spectra by using the action of a torus on the algebra $\Nil$.

Let $\mathcal{H}=(\C^{\ast})^{n}$, we define an $\mathcal{H}-$action on the algebra $\Nil$. Since the defining relations of the algebra $\Nil$ are homogeneous, we can define a torus $\mathcal{H}-$action on the algebra $\Nil$. Indeed, let $h=(\lambda_{1}, \lambda_{2}, \cdots, \lambda_{n})\in \mathcal{H}$, we can define the action of $h$ on the generators $e_{i}$ of $\Nil$ as follows:
\[
he_{i}=\lambda_{i}e_{i}
\]
for $i=1, 2, \cdots, n$. A non-zero element $x$ in the algebra $\Nil$ is called an $\mathcal{H}-$eigenvector of $\Nil$ if  we have $hx\in \C x$ for all $h\in\mathcal{H}$. 
An ideal $I$ of the algebra $\Nil$ is called $\mathcal{H}-$invariant if we have $hI =I$ for all $h\in \mathcal{H}$. We denote by $\mathcal{H}-Spec(\Nil)$ the set of all $\mathcal{H}-$invariant prime ideals of $\Nil$. First of all, by a theorem of Goodearl and Letzter about iterated skew polynomial rings, we know that this set is finite. 

In the one-parameter case, it was also proved by Gorelik in \cite{Gorelik} that the set $Spec_{\mathcal{H}}(U_{q}^{+}(\mathfrak g))$ is parameterized by the elements of the corresponding Weyl group. Thanks  to the fact that the algebra $\Nil$ is a cocycle twist of the algebra $U_{q}^{+}(\mathfrak g)$ and the cocycle twist establishes a one-to-one correspondence between the set of $\mathcal{H}-$invariant prime ideals of $U_{q}^{+}(\mathfrak g)$ and the set of $\mathcal{H}-$invariant prime ideals of $\Nil$, we have the following proposition.
\begin{prop}
The number of $\mathcal{H}-$invariant prime ideals of the algebra $\Nil$ is finite. In particular, the $\mathcal{H}-$invariant prime ideals of the algebra $\Nil$ are indexed by the elements of the Weyl group associated to $\mathfrak g$.
\end{prop}
\qed

Using Goodearl-Letzter stratification theory, the action of the torus $\mathcal{H}$ on the algebra $\Nil$ allows us to construct a partition of the prime spectrum $Spec(\Nil)$. Suppose that $J$ is an $\mathcal{H}-$invariant prime ideal of $\Nil$, we denote by $Spec_{J }(\Nil)$ the $\mathcal{H}-$stratum of the prime spectrum $Spec(\Nil)$ associated to the ideal $J$. Recall 
that we have $Spec_{J}(\Nil)\colon = \{P \in Spec(\Nil)\mid \cap_{h\in \mathcal{H}}hP = J\}$. And the $\mathcal{H}-$strata $Spec_{J}(\Nil)$ where $J\in \mathcal{H}-Spec(\Nil)$ form 
a partition of the prime spectrum $Spec(\Nil)$ as follows:
\[
Spec(\Nil) = \cup_{J\in \mathcal{H}-Spec(\Nil)} Spec_{J} (\Nil).
\]

Naturally, this partition induces a partition of the set $Prim(\Nil)$ of all (left) primitive ideals of $\Nil$. Namely, for any $J\in \mathcal{H}-Spec(\Nil)$, we shall denote primitive stratum corresponding to $J$ by  
\[
Prim_{J}(\Nil)\colon= Spec_{J}(\Nil)\cap Prim(\Nil).
\]
Then it is obvious to see that the $\mathcal{H}-$strata 
$Prim_{J} (\Nil)$ where $ J \in \mathcal{H}-Spec(\Nil)$
 form a partition of $Prim(\Nil)$ as follows:
\[
Prim(\Nil) = \cup_{J\in \mathcal{H}-Spec(\Nil)}Prim_{J} (\Nil).
\]

\subsection{Diximer-Moeglin equivalence, normal separation and catenarity}
In this subsection, we establish the fact that the algebra $\Nil$ satisfies the Diximer-Moeglin equivalence; and it has normal separation and is catenary. 

Let $k$ be any field and $A$ be a $k-$algebra. Recall that a prime ideal $P$ of a noetherian $k-$algebra $A$ is called rational provided the center of  Goldie quotient $Fract(A/P)$ is algebraic over the base field $k$.  Recall that a prime ideal $P$  is called a locally closed point in the prime spectrum $spec(A)$ under the Zariski topology, provided the singleton $\{P\}$ is a closed subset in some Zariski-neighborhood of the point $P$, i.e., the intersection of these primes properly containing $P$  is strictly larger than $P$.   We say that the algebra $A$ satisfies the {\bf Dixmier-Moeglin equivalence}  if the sets of rational prime ideals,  locally closed prime ideals, and primitive ideals all coincide. A useful criterion for Dixmier-Moeglin equivalence has been established in \cite{GLet} and we recall one of its versions, which works for many iterated skew polynomial rings with a torus action. 

Let 
\[
A = k[y_{1}][y_{2}; \tau_{2}; \delta_{2}] \cdots [y_{n} ;\tau_{n}; \delta_{n}]
\]
be an iterated skew polynomial ring over the field $k$.  Let $\mathcal{H}$ denote a group, whose elements act as $k-$algebra automorphisms on $A$ such that variables $y_{1}; \cdots ; y_{n}$ are all $\mathcal{H}-$eigenvectors.  For $1\leq i\leq  n$, let us set $A_{i} = k[y_{1}][y_{2}; \tau_{2}; \delta_{2}]\cdots [y_{i} ;\tau_{i}; \delta_{i}]$. The {\bf Theorem 4.7} in \cite{GLet} states that the algebra $A$ satisfies the Dixmier-Moeglin equivalence provided that the algebra $A$ satisfies the following three conditions:
\begin{enumerate}
\item There are infinitely many distinct eigenvalues for the action of $\mathcal{H}$ on $y_{1}$.
\item Each $\tau_{i}$  is a $k-$algebra automorphism of $A_{i-1}$  and each $\delta_{i}$ is a $k-$linear $\tau_{i}-$derivation of $A_{i-1}$.
\item For $2 \leq i  \leq n$, there exists $h_{i} \in  \mathcal{H}$  such that the restriction of $h_{i}$ to $A_{i-1}$ coincides with $\tau_{i}$ and the $h_{i}-$eigenvalue of $y_{i}$ is not a root of unity.
\end{enumerate}

It is obvious that the element $h$ acting on the algebra $\Nil$ as an algebra automorphism, which satisfies the conditions as outlined above. Thus via {\bf Theorem 4.4 \& Theorem 4.7} in \cite{GLet}, we can have the following result.
\begin{thm}
The algebra $\Nil$ satisfies the Dixmier-Moeglin equivalence. In particular, the primitive ideals of $A$ are precisely the prime ideals maximal within their $\mathcal{H}-$strata.
\end{thm}
\qed

Recall that the prime spectrum of an algebra $A$ is called  catenary
if for any two prime ideals $P\subsetneq Q \subset A$, all saturated chains of 
prime ideals from $P$  to $Q$ have the same length.  In addition, the prime spectrum of an algebra $A$ is said to have normal separation if the following condition holds: For any proper inclusion $P\subsetneq Q$ of prime ideals of $A$, the factor $Q/P$ contains a nonzero element which is normal in the factor algebra $A/P$.

Let $A$ be a group-graded algebra. According to \cite{Goodearl}, in order to show that the algebra $A$ has normal separation, it suffices to show that the algebra $A$ has graded-normal separation.  Recall that a graded-prime ideal $P$ of $A$ is any proper homogeneous ideal $P$ such that whenever $I, J$ are homogeneous ideals of $A$ with $IJ \subseteq P$, then either $I\subseteq P $ or $J\subseteq P$. We say that a graded algebra $A$ has graded-normal separation provided that for any  proper  inclusion $P \subsetneq  Q$ of  graded-prime  ideals of $A$,  there  exists a homogeneous element $c \in  Q-P$ which is normal modulo $P$. Once again, we can prove the following result using the fact that $U_{q}^{+}(\mathfrak g)$ has normal separation and $\Nil$ is a cocycle twist of $U_{q}^{+}(\mathfrak g)$.
\begin{prop}
The algebra $\Nil$ satisfies the graded normal separation, and therefore the normal separation. In particular, every non-zero prime ideal contains a non-zero normal element.
\end{prop}
\qed

We say that the Tauvel's  height  formula   holds in  an algebra   $A$  provided the following
\[
height(P) + GKdim(A/P) =  GKdim(A) 
\]
is true for any prime ideal $P\subsetneq A$.

To proceed, we need to recall a result from \cite{GLen}.
\begin{thm}
[{\bf Theorem  1.6.} in \cite{GLen}]  Let  R be  an affine,   noetherian,  Auslander-Gorenstein,    Cohen-Macaulay algebra  over  a field,   with finite  Gelfand-Kirillov   dimension.  If  $spec R$  
is normally  separated,   then  $R$  is  catenary.   If,  in  addition,  $R$  is  a  prime   ring,   then  Tauvel's  height formula   holds.
\end{thm}
\qed

Let $A$ be Noetherian $k-$algebra and let $M$  be a finitely generated module over $A$. Let us denote the Gelfand-Kirillov  and  homological dimensions  of $M$ by $GKdim(M)$, and respectively $hd(M)$. Let us denote the global homological dimension of $A$  by $gldim(A)$. Let us further denote by $injdim(A)$ the injective dimension of $A$. If $injdim(A)< \infty$, then $A$ is called Auslander-Gorenstein provided that $A$ satisfies the condition: for any integer $0\leq i\leq j$ and finitely generated (right) $A-$module $M$, we have $Ext^{i}(M,A)=0$ for all left $A-$submodules $N$ of $Ext^{j}(M,A)$.  If  $A$  is  an Auslander-Gorenstein  ring of finite global dimension, then $A$  is  called   Auslander-regular. Let us set $j(M)=min\{j:Ext^{j}(M, A) \neq 0\}$. The  ring $A$ is called Cohen-Macaulay (CM)  if $ j(M)+GKdim(M) =  GKdim(A)$ is true for  all  finitely generated $A-$modules $M$. 

Thanks to the iterated skew polynomial presentation of the algebra $\Nil$,  we can easily have the following result by using the {\bf Lemma } established in \cite{LS} repeatedly.
\begin{thm} The  two-parameter quantized  enveloping  algebra  $\Nil$ is an  affine   noetherian $\C-$algebra  with a finite  $GK-$dimension.  Moreover,  the algebra $\Nil$ is an  Auslander-regular, Cohen-Macaulay domain. 
\end{thm}
\qed

To tie the loose ends up, we have the following result.
\begin{thm} 
The  two-parameter quantized  enveloping  algebra  $\Nil$ is catenary. In particular, 
Tauvel's height formula holds for the algebra $\Nil$.
\end{thm}
\qed

\section{Prime ideals and Automorphism group of the algebra $\V$}

\subsection{The algebras $\U$ and $\V$}
Recall that $r, s$ are chosen from $\C^{\ast}$ such that $r,s$ are transcendental over $\Q$ and $r^{m}s^{n}=1$ implies that $m=n=0$. From reference \cite{BW1}, we shall recall the following construction of the algebra $\U$.
\begin{defn}
The two-parameter quantized enveloping algebra $\U$ is generated by the generators $e_{1}, e_{2}$ subject to the following relations:
\begin{eqnarray*}
e_{1}^{2}e_{2}-(r+s)e_{1}e_{2}e_{1}+rse_{2}e_{1}^{2}=0,\\
e_{1}e_{2}^{2}-(r+s)e_{1}e_{2}e_{1}+rse_{2}^{2}e_{1}=0.
\end{eqnarray*}
\end{defn}

Let us set the following new variables:
\begin{eqnarray*}
X_{1}=e_{1}, X_{2}=e_{2}, \\
X_{3}=X_{1}X_{2}-sX_{2}X_{1},\\
X_{3}^{\prime}=e_{1}e_{2}-r_{2}e_{2}e_{1}.
\end{eqnarray*}

Then we shall have the following identities:
\[
X_{1}X_{3}=rX_{3}X_{1},\,X_{2}X_{3}=r^{-1}X_{3}X_{2}.
\]

Let us further define some automorphisms $\tau_{1},\tau_{2},\tau_{3}$ and derivations
$\delta_{1},\delta_{2},\delta_{3}$ as follows: 
\begin{eqnarray*}
\tau_{2}(X_{1})=rX_{1},\delta_{1}(X_{1})=0\\
\tau_{3}(X_{1})=s^{-1}X_{1},\quad \tau_{3}(X_{3})=r^{-1}X_{3}\\
\delta_{3}(X_{1})=s^{-1}X_{3},\quad \delta_{3}(X_{3})=0.
\end{eqnarray*}

For the reader's convenience, we will state the following well-known result.
\begin{prop}
The algebra $\U$ can be presented as an iterated skew polynomial ring in that 
\[
\U \cong
\C[X_{1}][X_{3},\tau_{2},\delta_{2}][X_{2},\tau_{3},\delta_{3}].
\] 
In particular, the set $\{X_{1}^{i}X_{3}^{j}X_{2}^{k}|i, j, k\geq 0\}$ forms 
a PBW-basis of $\U$. The algebra $\U$ has a $GK-$dimension of $3$, and every 
prime ideal of $\U$ is completely prime.
\end{prop}
\qed 

It is easy to see that the algebra $\U$ is a special example of the downup algebras as defined 
in \cite{BT}. The primitives ideals of $\U$ were classified as a special case of the results in \cite{Praton}. The automorphism group of $\U$ was determined to be isomorphic to $(\C^{\ast})^{2}$ as a special case of the results recently established in \cite{CL}. In the reference \cite{T2}, using the deleting derivation algorithm \cite{Cauchon},  we have embedded the algebra $U^{+}_{r,s}(\mathfrak{s}\mathfrak{l}_{3})$ into a quantum torus, which enables us to determine all the derivations of $U^{+}_{r,s}(\mathfrak{s}\mathfrak{l}_{3})$ and compute the first Hochschild cohomology group for $U^{+}_{r,s}(\mathfrak{s}\mathfrak{l}_{3})$. Though most of the results on primitive ideals and prime ideals of for the algebra $\U$ are well-known, we will recall some of them in detail for the purpose of a smooth transition to the case of the algebra $\V$. To state the results for $\U$, we will follow the method used in reference \cite{M}. 

\subsection{Normal and central elements of the algebra $\U$}
  
First of all, we will recall some detailed information on the normal and central elements of the algebra $\U$. We shall verify that the normal elements of $\U$ are exactly scalar multiples of monomials in the variables $X_{3}=e_{1}e_{2}-re_{2}e_{1}$ and $X_{3}^{\prime}=e_{1}e_{2}-se_{2}e_{1}$.

Let $\alpha_{1},\alpha_{2}$ be the simple roots of $A_{2}$. We denote by $\mathcal{Q}$ the root lattice generated by $\alpha_{i}$ and denote by $\mathcal{Q}^{+}=\{\beta=\sum_{i=1}^{2}\beta_{i}\alpha_{i} 
\in \mathcal{Q}\mid \beta_{i}\geq 0\}$. For any $\beta \in \mathcal{Q}^{+}$, we define  
\[
\U_{\beta}=\C-span\{X_{1}^{l}X_{3}^{m}X_{2}^{n}\}
\]
with $\beta=l\alpha_{1}+m(\alpha_{1}+\alpha_{2})+n\alpha_{2}$.
\begin{prop}
$\U=\bigoplus_{\beta \in \mathcal{Q}} \U_{\beta}$. In particular, each
weight space $\U_{\beta}$ is finite dimensional.
\end{prop}
\qed

Note that the $\U$ also admits a filtration such that in the associated graded algebra, we have the following identities.
\begin{eqnarray*}
grX_{1}grX_{2}=sgrX_{2}grX_{1},\\
grX_{1}grX_{3}=rgrX_{3}grX_{1},\\
grX_{2}grX_{3}=r^{-1}grX_{3}grX_{2}.\\
\end{eqnarray*}
Let $X=X_{1}^{i}X_{3}^{j}X_{2}^{k}$, then we have 
\begin{eqnarray*}
grX_{1}grX=r^{j}s^{k}grXgrX_{1},\\
grX_{2}grX=r^{-j}s^{-i}grX grX_{2}.
\end{eqnarray*}

First of all, we have the following obvious lemma.
\begin{lem}
If $x$ is an $(r,s)-$central element, then $x$ is $(r,s)-$commuting with 
$X_{3}, X_{3}^{\prime}$.
\end{lem}
\qed

Now we have the following theorem on the normal elements and central elements of $\U$.
\begin{thm}
$N=\{c(X_{3})^{i}(X_{3}^{\prime})^{j}\mid c\in \C, i, j \geq 0\}$ is 
the set of all normal elements of $\U$. The center of $\U$ is $\C$.
\end{thm}
\qed

\subsection{Prime and primitive spectra of the algebra $\U$ and their stratifications}
The prime ideals of the quantum Heisenberg algebra (isomorphic to $U_{q}^{+}(sl_{3})$) 
were completely described in \cite{M}. The method used there can be carried over to
our algebra $\U$ with a minor modification. The difference is that the algebra $\U$ has a 
trivial center.  Using the method in \cite{M}, we will give a complete determination of all 
(completely) prime ideals of $\U$, and we will skip most of the details if no confusion arises. 

Note that the subalgebra of $\U$ generated by $X_{1},X_{3}$ or by
$X_{2}, X_{3}$ is a quantum plane. By abuse of notation, we will still
denote this subalgebra by $\C_{r,s}[X_{1}, X_{3}]$. For the prime ideals of
a quantum plane, we have the following results from \cite{M}.
\begin{lem}
[Lemma 2.1.2 in \cite{M}] Let $P$ be a prime ideal of $\C_{r,s}[X_{1},X_{3}]$ such that $P\cap
\C[X_{1}]=0$. Then $P$ is a principle ideal generated by the normal element $X_{3}$.
\end{lem}
\qed

\begin{prop}
[Proposition 2.1.3 in \cite{M}] All prime ideals of the quantum plane 
$\C_{r,s}[X_{1},X_{3}]$ are given as follows:
\begin{eqnarray*}
(0),\,(X_{1}),\,(X_{3});\\
(X_{1},X_{3});\\
(X_{1},X_{3}-\mu),\, (X_{3},X_{1}-\nu).
\end{eqnarray*}
for $\mu\nu\neq 0$.
\end{prop}
\qed 

Now we can derive some  results on prime ideals of the algebra $\U$ similar to those in \cite{M}.
\begin{lem}
Let $P$ be a prime ideal of $\U$ such that $P\cap \C_{r,s}[X_{1}, X_{3}]\neq
0$, then $P$ is the one of the following:
\begin{eqnarray*}
(X_{3}),\\
(X_{1}, X_{3}),\, (X_{2}, X_{3}),\\
(X_{1}, X_{2}-\mu, X_{3}),\, (X_{3}, X_{1}-\lambda, X_{2}),\\
(X_{3}, X_{1}, X_{2}).
\end{eqnarray*}
with $\lambda, \mu \in \C^{\ast}$.
\end{lem}
{\bf Proof:} (See the proof of {\bf Lemma 2.2.1} in \cite{M}) Let $P\subset \U$ be a prime ideal such that $P\cap
\C_{r,s}[X_{1}, X_{3}]\neq 0$. If $X_{1}\in P$, then $X_{3}\in P$. Suppose
$X_{3}\in P$, then the ideal $P/(X_{3})$ is a prime 
ideal of $\U/(X_{3})$ which is isomorphic to a quantum plane $\C[\overline{X_{1}}, \overline{X_{2}}]$. Thus we have that $P$ is one of the following prime ideals
of $\U$:
\begin{eqnarray*}
(X_{3}), (X_{1}, X_{3}),\,(X_{2}, X_{3}),\, (X_{3}, X_{1}-\lambda, X_{2}),\\
(X_{3}, X_{1}, X_{2}-\mu),\,(X_{1}, X_{2}, X_{3})
\end{eqnarray*}
with $\lambda\mu\neq 0$.

If $X_{1}-\lambda \in P$ for some $\lambda \neq 0$, then we have
$(X_{1}-\lambda)X_{2}-sX_{2}(X_{1}-\lambda)\in P$. Thus we have
$X_{3}+(s-1)\lambda X_{2}\in P$. Since $(1-s)\lambda \neq 0$, we have that 
$X_{2}\in P$. 

To finish the proof, we need to show that the situation $(X_{3}-\lambda, X_{1})\subset P$
can not happen for any $\lambda \neq 0$. Suppose we have that $(X_{3}-\mu, X_{1}) \subset P $ for some $\mu \neq 0$. Then we have $X_{3}\in P$, which implies that we shall
have $\mu=0$, a contradiction. Thus we have finished the proof. 
\qed 

Because of the symmetry between the variables $X_{3}, X_{3}^{\prime}$, we can also have the
following lemma.
\begin{lem}
Let $P$ be a prime ideal of $\U$ such that $P\cap \C_{r,s}[X_{1},
X_{3}^{\prime}]\neq 0$ and $X_{3}^{\prime}\in P$, then $P$ is 
given by replacing $X_{3}$ by $X_{3}^{\prime}$ in the previous 
lemma.
\end{lem}
\qed

So far, we have described all prime ideals of the algebra $\U$ which have
non-trivial intersections with the subalgebra $\C_{r,s}[X_{1}, X_{3}]$ or $\C_{r,s}[X_{1},X_{3}^{\prime}]$. Thus the problem is reduced to further describe those 
prime ideals of $\U$ which have trivial intersections with the subalgebras $\C_{r,s}[X_{1}, X_{3}]$ and $\C_{r,s}[X_{1}, X_{3}^{\prime}]$. We will show that such a prime ideal is indeed the zero ideal. The crucial point is to show that every non-zero two-sided ideal of $\U$ contains at least one nonzero normal element, and thus any non-zero prime ideal contains either $X_{3}$ or $X_{3}^{\prime}$. However, this is true due to the fact that $\U$ satisfies normal separation. Thus, we have the following proposition. 
\begin{prop}
Let $P\subset \U$ be a nonzero prime ideal of $\U$, then $P$ contains a non-zero normal element of $\U$. 
\end{prop}
\qed

Now we are ready to state the main result about prime ideals of the algebra $\U$.
\begin{thm}
All prime ideals of the algebra $\U$ are given as follows:
\begin{eqnarray*}
(0), (X_{3}),(X_{3}^{\prime}), \\
(X_{1}, X_{3}), (X_{2}, X_{3}), \\
(X_{3}, X_{1},X_{2}-\mu), (X_{3}, X_{1}, X_{2}), (X_{3},X_{1}-\lambda, X_{2}).
\end{eqnarray*}
\end{thm}
{\bf Proof:} Let $P\subset \U$ be a non-zero prime ideal of $\U$. Then
$P$ contains a non-zero normal element. Since $P$ is completely prime,
then $P$ contains either $X_{3}$ or $X_{3}$. Thus the theorem follows. 
\qed

To summarize everything, we state the following result.
\begin{thm}
The algebra $\U$ is catenary and satisfies the normal separation. In
addition, Tauvel's height formula holds for the algebra $\U$.
\end{thm}
\qed

By $Spec(\U)$, we denote the set of all prime ideals of $\U$. By $Prim(\U)$, we denote the set of all primitive ideals of $\U$. Note that the torus $\mathcal{H}=(\C^{\ast})^{2}$ is acting on the algebra $\U$ via algebra automorphisms. In order to describe the $\mathcal{H}-$stratification of the prime ideal spectrum, we need to single out all the $\mathcal{H}-$invariant prime ideals, which are indeed homogeneous ideals. It is easy to see that we have the following proposition.
\begin{prop}
The $\mathcal{H}-$invariant prime ideals of the algebra $\U$ are given as follows:
\[
(0), (X_{3}), (X_{3}^{\prime}), (X_{1}, X_{3}), (X_{2}, X_{3}), (X_{1}, X_{2}, X_{3}).
\]
\end{prop}
\qed

Now we describe the $\mathcal{H}-$stratification of the prime ideal spectrum for the algebra $\U$. Namely, we have the following theorem.
\begin{thm}
The $\mathcal{H}-$stratification of all the prime ideals of the algebra $\U$ is given as
follows:
\begin{eqnarray*}
Spec_{(0)}(\U)=\{(0)\};\\
Spec_{(X_{3})}(\U)=\{ (X_{3})\};\\
Spec_{(X_{3}^{\prime})}(\U)=\{(X_{3}^{\prime})\};\\
Spec_{(X_{1}, X_{3})}(\U)=\{(X_{1}, X_{3})\}\cup \{(X_{1}, X_{2}-\beta, X_{3})\mid \beta \in
\C^{\ast}\};\\
Spec_{(X_{2}, X_{3})}(\U) =\{(X_{2}, X_{3})\}\cup \{ (X_{1}-\alpha,X_{2}, X_{3})\mid \alpha \in \C^{\ast}\};\\
Spec_{(X_{1}, X_{2}, X_{3})}(\U) =\{ (X_{1}, X_{2}, X_{3})\}.
\end{eqnarray*}
\end{thm}
\qed 

From \cite{CL}, we have $G=Aut(\U)\cong (\C^{\ast})^{2}$. Note that the group $G$ is acting on the set of prime ideals. Therefore, we shall further have the following result.
\begin{prop}
The prime spectrum of the algebra $\U$ has $6$ $\mathcal{H}-$invariant strata which are
exactly $G-$orbits in the set $Spec(\U)$. 
\end{prop}  
 \qed

\subsection{The primitive spectrum of $\U$ and its $\mathcal{H}-$stratification}

Recall the primitive ideals of the downup algebras have been investigated by Praton in \cite{Praton}, where the primitive ideals are constructed in terms of the annihilators of 
irreducible representations of the downup algebras. In this section, we give a determination using the Dixmier-Moeglin equivalence. 
\begin{prop}
(See also \cite{Praton}) The primitive ideals of the algebra $\U$ are given as follows:
\begin{eqnarray*}
(0), (X_{3}), (X_{3}^{\prime}),\\
 (X_{3}, X_{1}, X_{2}-\beta),(X_{3}, X_{1}-\alpha, X_{2}), (X_{1},X_{2}, X_{3}).
\end{eqnarray*}
for $\alpha,\beta \in \C^{\ast}$.
\end{prop}
\qed

In addition, one can easily obtain the following result.
\begin{thm}
The $\mathcal{H}-$stratification of the primitive ideals is given as follows:
\begin{eqnarray*}
Prim_{(0)}(\U)=\{(0)\};\\
Prim_{(X_{3})}(\U)=\{(X_{3})\};\\
Prim_{(X_{3}^{\prime})}(\U)=\{(X_{3}^{\prime})\};\\
Prim_{(X_{1}, X_{3})}(\U)=\{(X_{1}, X_{2}-\beta, X_{3})\mid \beta \in \C^{\ast}\};\\
Prim_{(X_{2}, X_{3})}(\U)=\{(X_{1}-\alpha,X_{2}, X_{3})\mid \alpha \in \C^{\ast}\};\\
Prim_{(X_{1}, X_{2}, X_{3})}(\U) =\{ (X_{1}, X_{2}, X_{3})\}.
\end{eqnarray*}
\end{thm}
\qed

In particular, we shall further have the following result.
\begin{cor}
The primitive strata of $\U$ are exactly $G-$orbits in the set of all
primitive ideals. Furthermore, the ideal $(0)$ is a primitive ideal of $\U$, 
and thus the algebra $\U$ is a primitive ring.
\end{cor}
\qed

\subsection{Prime ideals and primitive ideals of the algebra $\V$}

First of all, we need to recall the following definition from reference \cite{BGH}.
\begin{defn}
The two-parameter quantized enveloping algebra $\V$ is defined to be the $\C-$algebra generated by the generators $e_{1}, e_{2}$ subject to the following relations:
\begin{eqnarray*}
e_{1}^{2}e_{2}-(r^{2}+s^{2})e_{1}e_{2}e_{1}+r^{2}s^{2}e_{2}e_{1}^{2}=0,\\
e_{1}e_{2}^{3}-(r^{2}+rs+s^{2})e_{2}e_{1}e_{2}^{2}
+rs(r^{2}+rs+s^{2})e_{2}^{2}e_{1}e_{2}-r^{3}s^{3}e_{2}^{3}e_{1}=0.
\end{eqnarray*}
\end{defn}

Now we further recall some basic properties of the algebra $\V$ as established in \cite{T3}. First of all, we need to fix more notation by setting the following new variables:
\begin{eqnarray*}
X_{1}=e_{1},\quad X_{2}=e_{3}=e_{1}e_{2}-r^{2}e_{2}e_{1},\\
X_{3}=e_{2}e_{3}-s^{-2}e_{3}e_{2},\quad
X_{4}=e_{2}.
\end{eqnarray*}

Now we recall the following result from reference \cite{T3}, whose proof is a straightforward calculation.
\begin{lem} The following identities hold.
\begin{enumerate}
\item $X_{1}X_{2}=s^{2}X_{2}X_{1}$;
\item $X_{1}X_{3}=r^{2}s^{2}X_{3}X_{1}$;
\item $X_{2}X_{3}=rsX_{3}X_{2}$;
\item $X_{1}X_{4}=r^{2}X_{4}X_{1}+X_{2}$;
\item $X_{2}X_{4}=s^{2}X_{4}X_{2}-s^{2}X_{3}$;
\item $X_{4}X_{3}=r^{-1}s^{-1}X_{3}X_{4}$.
\end{enumerate}
\end{lem}
\qed

In addition, let us define some algebra automorphisms $\tau_{2},\tau_{3}$, and $\tau_{4}$, and some derivations $\delta_{2},\delta_{3}$, and $\delta_{4}$ as follows: 
\begin{eqnarray*}
\tau_{2}(X_{1})=s^{-2}X_{1},\quad \delta_{2}(X_{2})=0,\\
\tau_{3}(X_{1})=r^{-2}s^{-2}X_{1},\,\tau_{3}(X_{2})=r^{-1}s^{-1}X_{2},\\
\delta_{3}(X_{1})=0,\,\delta_{3}(X_{2})=0,\\
\tau_{4}(X_{1})=r^{-2}X_{1},\,\tau_{4}(X_{2})=S^{-2}X_{2},\,\tau_{4}(X_{3})=r^{-1}s^{-1}X_{3},\\
\delta_{4}(X_{1})=-r^{-1}X_{2},\, \delta_{4}(X_{2})=X_{3},\, \delta_{4}(X_{3})=0.
\end{eqnarray*}
\qed

Then we have the following well-known result on a basis of the algebra $\V$.
\begin{thm}
The algebra $\V$ can be presented as an iterated skew polynomial ring. In
particular, we have the following result
\[
\U \cong
\C[X_{1}][X_{2},\tau_{2},\delta_{2}][X_{3},\tau_{3},\delta_{3}][X_{4}, \tau_{4}, \delta_{4}].
\] The set 
\[
\{X_{1}^{a}X_{2}^{b}X_{3}^{c}X_{4}^{d}|a, b, c, d \in \Z_{\geq 0}\}
\]
forms a PBW-basis of the algebra $\V$. In particular, the algebra $\V$ has a $GK-$dimension 
of $4$.
\end{thm}
\qed

Furthermore, using the graded algebra $gr(\V)$ of $\V$ associated to its obvious filtration, it is easy to verify directly the following result.
\begin{cor}
The center of the algebra $\V$ is reduced to the base field $\C$.
\end{cor}
\qed

Following the convention in the case of $U_{q}^{+}(B_{2})$, we may also denote the element $X_{3}$ by $Z$. In addition, we shall further set the following variable:
\[
Z^{\prime}=(X_{1}(X_{3}+(s^{-2}-r^{-1}s^{-1})X_{2}X_{4})-s^{4}(X_{3}+(s^{-2}-r^{-1}s^{-1})X_{2}X_{4})X_{1}).
\]

For convenience, let us further set a new variable $W=X_{3}+(s^{-2}-r^{-1}s^{-1})X_{2}X_{4}$. Then we shall have the following lemma, which can be verified by brutal force.

\begin{lem} The following identities can be verified to hold in $\V$.
\begin{enumerate}
\item $X_{1}W=r^{2}s^{2}WX_{1}+(1-r^{-1}s)X_{2}^{2}$;\\
\item $X_{2}W=s^{2}WX_{2}$;\\
\item $X_{3}W=WX_{3}$;\\
\item $X_{4}W=s^{-2}WX_{4}$;\\
\item $X_{1}Z^{\prime}=r^{2}s^{2}Z^{\prime}X_{1}$;\\
\item $X_{2}Z^{\prime}=Z^{\prime}X_{2}$;\\
\item $X_{3}Z^{\prime}=r^{-2}s^{-2}Z^{\prime}X_{3}$;\\
\item $X_{4}Z^{\prime}=r^{-2}s^{-2}Z^{\prime}X_{4}.$
\end{enumerate}
\end{lem}
\qed

As a result of the previous identities, we know that both the elements $Z$ and $Z^{\prime}$ are $r-s-$central elements of the algebra $\V$. Thus both $Z$ and $Z^{\prime}$ are normal elements of the algebra $\V$ because normal elements of the algebra $\V$ are just $r-s-$central elements of the algebra $\V$. In addition, we shall have the following lemma, which describes all the normal elements of the algebra $\V$.
\begin{lem}
The normal elements of $\V$ are of the form $\alpha Z^{m}(Z^{\prime})^{n}$ where $\alpha \in \C$ and $m, n \in \Z_{\geq 0}$.
\end{lem}
{\bf Proof:} Let $u$ be a normal element of $\V$. Then $u$ is an $r-s$central element of $\V$. Let us embed the algebra $\V$ into the subalgebra $\C_{r,s}[e_{1}^{\pm 1}, e_{3}^{\pm 1}, Z, Z^{\prime}]$ of the Goldie quotient ring of $\V$. Therefore, the element $u$ takes the following format
\[
u=\sum \alpha_{a,b,m, n} e_{1}^{a}e_{3}^{b}Z^{m}(Z^{\prime})^{n}.
\]
Since the element $u$ is $r-s-$commuting with the elements $e_{1}, e_{2}$, and $e_{3}=e_{1}e_{2}-r^{2}e_{2}e_{1}$, the element $u$ shall $r-s-$commute 
with the elements $e_{1}$ and $e_{3}$. Via direct calculations, we can prove that
\[
u=\alpha e_{1}^{a}e_{3}^{b}Z^{m}(Z^{\prime})^{n}.
\]
Furthermore, using the fact that the element $u$ is $r-s-$commuting with $e_{2}$, via direct calculations, we can prove that the element $u$ indeed takes the following format
\[
u=\alpha Z^{m}(Z^{\prime})^{n}
\]
for some $\alpha \in \C$ as desired. Therefore, we have finished the proof.
\qed
 
Furthermore, using the general properties previously established for $U^{+}_{r,s}(\mathfrak g)$, we shall have the following result.
\begin{prop}
Every non-zero prime ideal of $\V$ either contains $Z$ or $Z^{\prime}$. In particular, the algebra $\V$ satisfies the normal separation, is catenary; and the Tauvel's Formula and the Dixmier-Moeglin equivalence hold. 
\end{prop}
\qed

The next proposition relates the algebra $\V$ to the algebra $\U$ as previously discussed.
\begin{prop}
Let $(Z)$ be the two-sided ideal of the algebra $\V$ generated by the normal element $Z$, then we have the following
\[
\V/(Z)\cong \U
\]
as an algebra.
\end{prop}
\qed

In addition, using the theory developed in \cite{Goodearl, GLet}, we can easily verify the following proposition.
\begin{prop}
The ideal $(Z^{\prime})$ generated by the normal element $Z^{\prime}$ is a prime ideal of the algebra $\V$. Furthermore, the ideal $(Z^{\prime})$ is a primitive ideal of the algebra $\V$.
\end{prop}
{\bf Proof:} Since the algebra $\V$ can be embedded into a subalgebra $\C_{r,s}[e_{1}^{\pm 1}, e_{3}^{\pm 1}, Z, Z^{\prime}]$ of the Goldie quotient ring of $\V$, the quotient algebra $\V/(Z^{\prime})$ can be embedded into the ring $\C_{r,s}[e_{1}^{\pm 1}, e_{3}^{\pm 1}, Z]$ which is obviously a domain. Therefore, we have proved that $(Z^{\prime})$ is a (completely) prime ideal of $\V$. In addition, we can verify that the center of the localization of $\V/(Z^{\prime})$ with respect to the Ore set associated to the prime ideal $Z^{\prime}$ (see \cite{Goodearl} for the information on this Ore set) is the base field $\C$. Using {\bf Theorem 5.3.} in \cite{Goodearl}, we can prove that $(Z^{\prime})$ is the only prime ideal in the stratum associated to the $\mathcal{H}-$invariant ideal $(Z^{\prime})$. As a result, the prime ideal $(Z^{\prime})$ is also a primitive ideal of $\V$.
\qed

In addition, due to the fact that the algebra $\V$ satisfies the normal separation and $\alpha Z^{m}Z^{\prime n}$ are the only normal elements of $\V$, we have the following result.
\begin{prop}
The ideal $(0)$ is the only element in $Spec_{(0)}(\V)$. Therefore, we know that that the prime ideal $(0)$ is a primitive ideal of the algebra $\V$. Therefore, the algebra $\V$ is primitive.
\end{prop}
\qed

In addition, it is easy to see that we have the following.
\begin{thm}
\begin{enumerate}
\item The algebra $\V$ has $8$ $\mathcal{H}-$invariant prime ideals, which are listed as follows:
\begin{eqnarray*}
(0), (Z), (Z^{\prime}), (e_{3}), 
\\(\overline{e_{3}}), (e_{1}), (e_{2}), (e_{1}, e_{2}).
\end{eqnarray*}
\item The stratification of the prime ideals of the algebra $\V$ is given as follows:
\begin{eqnarray*}
Spec_{(0)}(\V)=\{(0)\};\\
Spec_{(Z)}(\V)=\{(Z)\};\\
Spec_{(Z^{\prime})}(\V)=\{(Z^{\prime})\};\\
Spec_{(e_{3})}(\V)=\{(e_{3})\};\\
Spec_{(\overline{e_{3}})}(\V)=\{(\overline{e_{3}})\};\\
Spec_{(e_{1})}(\V)=\{(e_{1})\}\cup \{(e_{1}, e_{2}-\mu)\mid \mu \in \C^{\ast}\};\\
Spec_{(e_{2})}(\V)=\{(e_{2})\}\cup \{(e_{1}-\lambda, e_{2})\mid \lambda \in \C^{\ast}\};\\
Spec_{(e_{1}, e_{2})}(\V)=\{(e_{1}, e_{2})\}.
\end{eqnarray*}
\end{enumerate}
\end{thm}
\qed

As a result, we have the following description of the primitive ideals of the algebra $\V$.
\begin{thm}
The stratification of the primitive ideals of the algebra $\V$ is given as follows:
\begin{eqnarray*}
Prim_{(0)}(\V)=\{(0)\};\\
Prim_{(Z)}(\V)=\{(Z)\};\\
Prim_{(Z^{\prime})}(\V)=\{(Z^{\prime})\};\\
Prim_{(e_{3})}(\V)=\{(e_{3})\};\\
Prim_{(\overline{e_{3}})}(\V)=\{(\overline{e_{3}})\};\\
Prim_{(e_{1})}(\V)=\{(e_{1}, e_{2}-\mu)\mid \mu \in \C^{\ast}\};\\
Prim_{(e_{2})}(\V)=\{(e_{1}-\lambda, e_{2})\mid \lambda \in \C^{\ast}\};\\
Prim_{(e_{1}, e_{2})}(\V)=\{(e_{1}, e_{2})\}.
\end{eqnarray*}
\end{thm}
\qed

\subsection{The automorphism group of the algebra $\V$}

As the defining relations of the algebra $\V$ are homogeneous in the 
given generators $e_{1}, e_{2}$, there is an $\N-$grading
on the algebra $\V$ obtained by assigning to $e_{i}$ the degree 1. Let
\[
\V=\bigoplus_{i\in \N} \V_{i}
\]
be the corresponding decomposition, with $\V_{i}$ being 
the subspace of homogeneous elements of degree i. 
In particular, we have that $\V_{0} = \C$ and $\V_{1}$ is the two-dimensional space spanned by 
the generators $e_{1}, e_{2}$. For any $t \in \N$, we can set $\V_{\geq t} = \bigoplus_{i\geq t}\V_{i}$ and define $V_{\leq t}=\bigoplus _{i\leq t}\V_{i}$. We say that the nonzero element $u\in
\V_{\geq t}-\V_{\leq t-1}$ has degree $t$, and write $deg(u) = t$. Since
the algebra $\V$  is a domain, we have $deg(uv) = deg(u) + deg(v)$ for $u, v \neq
0$.

Now we recall the following definition from reference \cite{LL}.
\begin{defn}
Let $A = \bigoplus_{i \in \N} A_{i}$ be an $\N-$graded $\C-$algebra with $A_{0} = \C$,  which is generated as 
an algebra by $A_{1} = \C x_{1} \oplus \cdots \oplus \C x_{n}$. 
If for each $i \in \{1, . . . , n\}$, there exists 
$0 \neq a \in A$ and a scalar $q_{i}$ such that $x_{i}a =
q_{i}ax_{i}$, then we say that $A$ is an
$\N-$graded algebra with enough $q-$commutation relations.
\end{defn}

According to \cite{LL}, it is easy to see that we have the following results.
\begin{prop}
The algebra $\V$ endowed with the gradation as just defined, is a 
connected $\N-$graded algebra with enough $q-$commutation relations.
\end{prop}
\qed 

\begin{cor}
Let $\sigma \in Aut(\V)$ and $x\in \V_{d}-0$, then
$\sigma(x)=y_{d}+y_{y>d}$ where $y_{d}\in \V_{d}-0$ and 
$y_{d}\in \V_{\geq d+1}$.
\end{cor}
\qed

In order to proceed, we need to establish some basic identities. First of all, recall that we have denoted  by $\overline{e_{3}}=e_{1}e_{2}-s^{2}e_{2}e_{1}$. We have the following lemma.
\begin{lem} Let us assign degree $1$ to both of the generators $e_{1}, e_{2}$. Then we have 
\begin{eqnarray*}
Z^{\prime}&=&r^{2}s^{2}(1-r^{-2}s^{2})Ze_{1}+(1-r^{-1}s)(1-r^{-2}s^{2})e_{3}e_{1}e_{2}\\
&&+r^{-2}s^{2}(1-r^{-1}s)e_{3}^{2}\\
&=&r^{2}(1-r^{-2}s^{2})(s^{2}Z+(1-r^{-1}s)e_{3}e_{2})e_{1}+(1-r^{-1}s)e_{3}^{2}\\
&=& (rs)^{-1}(1-r^{-1}s)e_{3}\overline{e_{3}}+rs(1+r^{-1}s)Ze_{1}\\
&=& (rs)^{-1}(1-r^{-1}s)(\overline{e_{3}})^{2}+ue_{1}
\end{eqnarray*} 
for some homogeneous element $u\in \V$ of degree $3$.
\end{lem}
{\bf Proof:} These identities can be verified via brutal force calculations, and we will not repeat them here.
\qed

\begin{lem}
Let $\sigma\in Aut(\V)$ be an algebra automorphism of $\V$, then we have the following
\[
\sigma(Z)=\lambda Z,\quad \sigma(Z^{\prime})=\mu Z^{\prime}
\]
for some $\lambda, \mu \in \C^{\ast}$.
\end{lem}
{\bf Proof:} Since the normal elements of the algebra $\V$ are just $\alpha Z^{m}(Z^{\prime})^{n}$ and the normal elements are sent to normal elements by the algebra automorphism $\sigma$, we have that $\sigma(Z)=\lambda Z$ or $\sigma(Z)=\lambda  Z^{\prime}$ for some $\lambda \in \C^{\ast}$. Via the commuting relations between $e_{1}, e_{2}$ and $Z, Z^{\prime}$, we can further prove that $\sigma(Z)=\lambda Z$ and $\sigma (Z^{\prime})=\mu Z^{\prime}$ as desired.
\qed

\begin{thm}
Let $\sigma\in Aut(\V)$. Then we have the following
\[
\sigma(e_{1})=\alpha_{1} e_{1}, \quad \sigma(e_{2})=\alpha_{2}e_{2}
\]
for some $\alpha_{1}, \alpha_{2}\in \C^{\ast}$. In particular, we have the following
\[
Aut(\V)\cong (\C^{\ast})^{2}.
\]
\end{thm}
{\bf Proof:}(Our proof follows the argument used in \cite{L1} for the case of $U_{q}^{+}(B_{2})$) Let $\sigma\in Aut(\V)$ be an algebra automorphism of the algebra $\V$. Based on the graded arguments, we shall have the following
\[
\sigma(X_{i})=\alpha_{i}X_{i}+u_{i}
\]
for $i=1, 2, 3, 4$ and the $u_{i}$ is of higher degree than the $X_{i}$ unless $u_{i}$ is zero. 

For $i=1, 2$ and $3$, let us set $d_{i}=deg(\sigma(e_{i}))$ and $d_{4}=\sigma(s^{2}Z+(1-r^{-1}s)e_{3}e_{2})$. Let us further set $d_{3}^{\prime}=deg(\sigma(\overline{e_{3}}))$. We have that $d_{1}, d_{2}\geq 1$ and $d_{4}\geq 3$ and $d_{3}, d_{3}^{\prime}\geq 2$. It suffices to show that $d_{1}=d_{2}=1$.

Since $Z$ is fixed by $\sigma$ up to a non-zero scalar, we have that $deg(\sigma(Z))=3$. Suppose that we have $d_{2}+d_{3}>3$. Since $Z^{\prime}$ is also fixed by $\sigma$ up to a 
non-zero scalar, we have that $d_{1}+d_{4}=d_{1}+d_{2}+d_{3}=2d_{3}>4$, which implies $d_{1}+d_{2}=d_{3}>2$. Thus we have that 
$d_{3}+d_{3}^{\prime}>4$.  So we have that $3+d_{1}=d_{3}+d_{3}^{\prime}$, which implies that $d_{2}+d_{3}^{\prime}=3$. This implies 
that $d_{2}=1$ and $d_{3}^{\prime}=2$. If $d_{2}+d_{3}=3$, then we have $d_{2}=1, d_{3}=2$. So we always have $d_{2}=1$ and either $d_{3}=2$ or 
$d_{3}^{\prime}=2$. In either case, we shall have $d_{1}=1$ as desired. Using the commuting relations between $(e_{1}), e_{2}$ and $Z$, and 
the fact that $\sigma$ fixes $Z$ up to a non-zero scalar, we can finish the proof of the theorem. 
\qed

It is obviously that the automorphism group $Aut(\V)$ acts on the set of primitive ideals. In particular, we have the following result.
\begin{prop}
The $Aut(\V)$-orbits within the primitive spectrum $Prim(\V)$  coincide with the 
$\mathcal{H}-$strata of $Prim(\V)$. That is, the action of $Aut(\V)$ on $Prim(\V)$ 
has exactly 8 orbits
\begin{eqnarray*}
Prim_{(0)}(\V)=\{(0)\};\\
Prim_{(Z)}(\V)=\{(Z)\};\\
Prim_{(Z^{\prime})}(\V)=\{(Z^{\prime})\};\\
Prim_{(e_{3})}(\V)=\{(e_{3})\};\\
Prim_{(\overline{e_{3}})}(\V)=\{(\overline{e_{3}})\}\\
Prim_{(e_{1})}(\V)=\{(e_{1}, e_{2}-\mu)\mid \mu \in \C^{\ast}\};\\
Prim_{(e_{2})}(\V)=\{(e_{1}-\lambda, e_{2})\mid \lambda \in \C^{\ast}\};\\
Prim_{(e_{1}, e_{2})}(\V)=\{(e_{1}, e_{2})\}.
\end{eqnarray*}
\end{prop}
\qed

\end{document}